\documentclass[titlepage,11pt]{article}
\usepackage{latexsym}
\usepackage{amsfonts}
\usepackage{amsmath}
\newcommand\blackslug{\hbox{\hskip 1pt \vrule width 4pt height 8pt depth 1.5pt
        \hskip 1pt}}
\newcommand\bbox{\hfill \quad \blackslug \bigbreak}

\def\LL{,\ldots,}

\title{A note on the Gy\'arf\'as-Sumner conjecture}
\author{Tung Nguyen\thanks{Supported by AFOSR grant  FA9550-22-1-0234 and by NSF grant  DMS-2154169.}\\
Princeton University,\\ Princeton, NJ 08544, USA
\and
Alex Scott\thanks{Research supported by EPSRC grant EP/X013642/1.}\\
Mathematical Institute,\\ University of Oxford,
\\
Oxford OX2 6GG, UK
\and
Paul Seymour\footnotemark[1]\\
Princeton University,\\ Princeton, NJ 08544, USA}

\date{January 17, 2023; revised \today}

\newtheorem{thm}{}[section]

\newcommand{\Proof}{\noindent{\bf Proof.}\ \ }

\begin{document}
\maketitle
\begin{abstract}

The Gy\'arf\'as-Sumner conjecture says that for every tree $T$ and every integer $t\ge 1$, if $G$ is a graph with no clique of size $t$
and with sufficiently large chromatic number, then $G$ contains an induced subgraph isomorphic to $T$. This remains open, but we prove that
under the same hypotheses, $G$ contains a subgraph $H$ isomorphic to $T$ that is ``path-induced''; that is, for some distinguished vertex~$r$, every path of $H$ with one end $r$ is an induced path of $G$.

\end{abstract}

\section{Introduction}

The Gy\'arf\'as-Sumner conjecture says~\cite{gyarfastree,sumner}:
\begin{thm}\label{GSconj}
{\bf Conjecture: }For every tree 
$T$, and every integer $t\ge 1$, if $G$ is a graph with 
no induced subgraph isomorphic to $T$, and with no 
clique of size $t$,
then its chromatic number is bounded.
\end{thm}
This has been proved for a few families of trees (for instance \cite{distantstars, gyarfasprob, gst, kierstead, kierstead2, scott, newbrooms}, and 
see~\cite{chibounded} for a survey), but remains open in general.

Let $T$ be a tree and let $r$ be a vertex of $T$ (we call $(T,r)$ a {\em rooted tree}).  Suppose that $\phi$ is an isomorphism from $T$ to a subgraph of a graph $G$ (not necessarily induced).
We say that $\phi$ is a {\em path-induced copy of $(T,r)$} in $G$ if for every path $P$ of $T$ with one end $r$, the image of $P$ under $\phi$ is
an induced path in $G$.
The notion of being path-induced lies
partway between being a subgraph and being an induced subgraph.

\ref{GSconj} remains out of reach.  However, we will show here that it is true if instead of excluding an {\em induced} tree we
exclude
a {\em path-induced} tree.  Our main result is the following.
%
\begin{thm}\label{mainthm}
For every rooted tree $(T,r)$, and every integer $t\ge 1$, if $G$ is a graph with
no path-induced copy of $(T,r)$, and with no
clique of size $t$,
then its chromatic number is bounded.
\end{thm}

We prove \ref{mainthm} in the next section, and give further discussion in the final section.

\section{The proof}

We denote the chromatic number of $G$ by $\chi(G)$, and for $X\subseteq V(G)$, we write $\chi(X)$ for $\chi(G[X])$.
We need the following lemma. It is a special case of a theorem of~\cite{distantstars}, but we give the proof since it is short.
For $c\ge 1$, we say $X\subseteq V(G)$ is a {\em $c$-creature} of $G$ if every vertex in $X$ has fewer than $c$ neighbours in $V(G)\setminus X$.

\begin{thm}\label{creature}
For all integers $a,c\ge 0$, if $G$ is a graph with $\chi(G)>ac$, and $X\subseteq V(G)$ is a $c$-creature with $\chi(X)\le a$, then
$\chi(G\setminus X)=\chi(G)$.
\end{thm}
\Proof Suppose that $\chi(G\setminus X)<\chi(G)$, and let $\kappa:V(G)\setminus X\rightarrow \{1\LL \chi(G)-1\}$
be a colouring of $G\setminus X$ with $\chi(G)-1$ colours.
Since $\chi(X)\le a$, there is a partition $X_1\LL X_a$ of $X$ into $a$ stable sets. For $1\le i\le a$, let 
$$J_i=\{(i-1)c+1,(i-1)c+2\LL ic\}.$$ 
For each $v\in X_i$, choose $\kappa(v)\in J_i$ different from $\kappa(u)$
for each neighbour $u\in V(G)\setminus X$ of $v$ (this is possible since $|J_i|=c$ and $v$ has fewer than $c$ neighbours in $V(G)\setminus X$). 
Thus we have extended $\kappa$ to a $(\chi(G)-1)$-colouring of $G$, a contradiction. 
This proves \ref{creature}.~\bbox

If $\phi$ is a path-induced copy of some $(T,r)$, we denote by $V(\phi)$ the vertex set of its image, that is, 
the set of all vertices $\phi(h)\;(h\in V(T))$.
For integers $d\ge 2$ and $k\ge 1$, let $(T^k_d, r)$ be the rooted tree in which the root $r$ has degree $d$, 
every vertex has degree $d$ or $1$, and every path from 
$r$ to a leaf has length $k$. We will prove:

\begin{thm}\label{getcreature}
Let $k\ge 1$, $d\ge 2$, and $\tau\ge 0$ be integers. Then there exists $K$ with the following property.
Let $G$ be a graph, such that for every $v\in V(G)$,
$\chi(N(v))\le \tau$, where $N(v)$ is the set of neighbours of $v$. Let $v\in V(G)$. Then either $G$ admits a 
path-induced copy $\phi$ of $(T_d^k,r)$ with $\phi(r) = v$, or there
exists a $(1+d+d^2+\cdots+d^{k-1})$-creature  $X$ of $G$ with $v\in X$ such that $\chi(X)\le K$.
\end{thm}
\Proof
For each $j\ge 1$, define $c(j) = 1+d+d^2+\cdots+d^{j-1}$. Define $f(1) = 1$, and inductively for $j\ge 2$,
define $f(j) = f(j-1)c(j-1)+\tau$. We will prove by induction on $k$ that $K=f(k)$ satisfies the theorem.

If $k=1$, and $v$ has $d$ neighbours in $G$, then there is a path-induced copy $\phi$ of $(T_d^k,r)$
with $\phi(r) = v$; and otherwise $\{v\}$ is a $d$-creature. Thus we may assume that $k\ge 2$ and the result holds for $k-1$.
Let $M(v)=V(G)\setminus (N(v)\cup \{v\})$. 
Choose $A\subseteq N(v)$ maximal such that for each $a\in A$,
there is a path-induced copy $\phi_a$ of $(T_d^{k-1},r)$ such that 
\begin{itemize}
\item $\phi_a(r) = a$ and $V(\phi_a)\subseteq \{a\}\cup M(v)$, for each $a\in A$; and
\item the sets $V(\phi_a)\; (a\in A)$ are pairwise disjoint.
\end{itemize}
If $|A|\ge d$ then $G$ admits a 
path-induced copy $\phi$ of $(T_d^k,r)$ with $\phi(r) = v$, as required, so we may assume that $|A|<d$. Consequently
the union of the sets $V(\phi_a) \; (a\in A)$ has cardinality at most $(d-1)(1+d+d^2+\cdots d^{k-1})=d^k-1$.
Let us denote this union by $W$.

For each $u\in N(v)\setminus A$, from the inductive hypothesis, applied to the subgraph induced on $(\{u\}\cup M(v))\setminus W$,
there is a $c(k-1)$-creature $X_u$ of this subgraph with $u\in X_u$ and with $\chi(X_u)\le f(k-1)$. Let $X$ be the union of $\{v\}$
and all the sets $X_u\; (u\in N(v)\setminus A)$. We claim that $X$ satisfies the theorem. 
\\
\\
(1) {\em $X$ is a $c(k)$-creature.}
\\
\\
First, since $X$ contains all vertices of $N(v)\setminus A$, it follows that 
$v$ has $|A|<d\le c(k)$ neighbours in $V(G)\setminus X$.
Every other vertex $x\in X$ belongs to one of the sets $X_u$ where $u\in N(v)\setminus A$, and so has fewer than 
$c(k-1)$ neighbours in $M(v)\setminus (X_u\cup W)$, and consequently has fewer than $c(k-1)$ in 
$V(G)\setminus (X\cup W)$. Moreover, it has at most $d^k$ neighbours in $W$ because $|W|\le d^k$, and so has in total
fewer than $c(k-1)+d^k=c(k)$ neighbours in $V(G)\setminus X$. This proves (1).
\\
\\
(2) {\em $\chi(X)\le f(k)$.}
\\
\\
For each $B\subseteq N(v)\setminus A$, let $X_B$ be the union of the sets $(X_b\setminus \{b\})\;(b\in B)$. Thus 
$X_B\subseteq M(v)\setminus W$. Suppose first that $\chi(X_{N(v)\setminus A})> f(k-1)c(k-1)$, and choose 
$B\subseteq N(v)\setminus A$ minimal such that $\chi(X_B))>  f(k-1)c(k-1)$. Since $B\ne \emptyset$, there exists 
$b\in B$; but $X_b\setminus \{b\}$ 
is a $c(k-1)$-creature of $G[X_B]$, with chromatic number at most $f(k-1)$, contrary to \ref{creature} and the minimality
of $B$. This proves that $\chi(X_{N(v)\setminus A)})\le  f(k-1)c(k-1)$. Since $\chi(N(v))\le \tau$, it follows that
$\chi(X)\le  f(k-1)c(k-1)+\tau=f(k)$. This proves (2).

\bigskip

From (1) and (2), $X$ satisfies the theorem. This proves \ref{getcreature}.~\bbox

Now we can deduce \ref{mainthm}, which we restate:
\begin{thm}\label{mainthm2}
For every rooted tree $(T,r)$ and every integer $t\ge 1$, if $G$ is a graph with no clique of size $t$
and with sufficiently large chromatic number, then $G$ admits a path-induced copy of $(T,r)$.
\end{thm}
\Proof
We may assume that $(T,r)$ equals $(T_d^k,r)$ for some choice of $d\ge 2$ and $k\ge 1$. We proceed by induction on $t$;
and so may assume that there exists $\tau\ge 0$ such that if $G$ is a graph with no clique of size $t-1$,
and $\chi(G)>\tau$, then $G$ admits a path-induced copy of $(T,r)$.
Choose $f(k)$ as in \ref{getcreature}.
We claim that if $G$ is a graph with no clique of size $t$ that does not admit a 
path-induced copy of $(T,r)$ then $\chi(G)\le (1+d+d^2+\cdots+d^{k-1})f(k)$. We may assume that
$G$ is critical with its chromatic number; that is, for every nonempty subset $X\subseteq V(G)$, $\chi(G\setminus X)<\chi(G)$.
By \ref{getcreature}, 
there is a $(1+d+d^2+\cdots+d^{k-1})$-creature  $X$ of $G$ with $v\in X$ such that $\chi(X)\le f(k)$. 
Since $G$ has no clique of size $t$, for every vertex $v$, $G[N(v)]$ has no clique of size $t-1$, and so $\chi(G[N(v)])\le \tau$.
Since $\chi(G\setminus X)<\chi(G)$, \ref{creature} implies that 
$$\chi(G)\le (1+d+d^2+\cdots+d^{k-1})f(k).$$ 
This 
proves \ref{mainthm2}.~\bbox

\section{Strengthenings}

In this final section, we discuss possible strengthenings of \ref{mainthm}.

One can refine \ref{mainthm} a little, and we first give a sketch of this. Start with a graph $G$ with bounded clique number and large
(really really huge!) chromatic number.
For fixed $k,d$, we can apply \ref{mainthm} to get a path-induced copy of $(T_D^k,r)$ in $G$  for some huge value of $D$,
and then use Ramsey arguments to get a path-induced copy of $T_d^k$ where we have some control over the edges that stop this subgraph
being induced. For instance, since there is a bound on the maximum size of a clique, each vertex of the tree that has children has
a large (say size $D'$) set of children that is stable in $G$; and the set of vertices such that all their ancestors belong to 
the selected stable subsets forms a path-induced copy of $(T_{D'}^k,r)$ in which for every vertex, its set of children is stable in $G$.
We can also
arrange, using the bipartite Ramsey theorem, that for every two vertices of the tree at the same height, the children of the first are
either completely adjacent, or completely nonadjacent, to the children of the second. And then 
we can get a path-induced copy of $(T_{D''}^k,r)$, for some huge $D''$, such that for each 
vertex, its set of grandchildren is stable in $G$. And so on: we can arrange that for each $i$,
the set of vertices at distance $i$ from the root is a stable set. (Let us call this being {\em level-stable}.)

We can also
arrange, using the bipartite Ramsey theorem, that for every two vertices of the tree that are not leaves, even if their height is different, the children of the first are
either completely adjacent, or completely nonadjacent, to the children of the second. (The argument here is more tricky: it is important to fix up pairs in the right order, 
but we omit the details.) 

But we can go further. Say two vertices $u,v$ of the tree are {\em incomparable} if neither is an ancestor of the other; and 
let $d(u,v)$ denote the distance between $u,v$ in the tree. If $u,v$ belong to the tree, let $w$ be their 
``join'' (their common ancestor furthest from the root), and let $a=d(u,r), b=d(v,r)$ and $c=d(w,r)$. The triple $(a,b,c)$
describes the pair $u,v$ up to isomorphisms of the tree. But we need a little more information. For each vertex, choose a linear order
of its set of children. So now, if $u,v$ are incomparable, then they descend from different children (say $u', v'$ respectively) 
of their join $w$, and one of these is earlier than the other in the linear order of the children of $w$. If $u'$ is earlier than $v'$
we say $u$ is earlier than $v$. Let us say that the {\em type} of an unordered pair $\{u,v\}$ (where $u,v$ are incomparable) 
is the triple
$(a,b,c)$ defined as before, where $u$ is earlier than $v$. 

Let us say a path-induced copy of $T_d^r$
in $G$ is {\em type-uniform} if the adjacency of each incomparable pair of vertices depends only on their type; in other words, if $\{u,v\}$
and $\{u',v'\}$ have the same type, then they are both adjacent or both nonadjacent pairs.
One can use more Ramsey arguments (we omit the details, which are straightforward) to arrange,
again by reducing $D$, that the adjacency of each pair of vertices depends only on their type; in other words, if $\{u,v\}$
and $\{u',v'\}$ have the same type, then they are both adjacent or both nonadjacent pairs. 
In conclusion, then, we deduce:
\begin{thm}\label{type-uniform}
For all $k,t\ge 1$ and $d\ge 2$, if $G$ is a graph with no clique of size $t$
and with sufficiently large chromatic number, then $G$ admits a path-induced, level-stable, type-uniform copy of $(T_d^k,r)$.
\end{thm}

\medskip

Here is another way in which we might strengthen \ref{mainthm}:
can we obtain polynomial bounds?  There is an analogous problem where we exclude the complete bipartite graph $K_{t,t}$ as a subgraph, 
instead of excluding $K_t$, and 
the following was shown in~\cite{poly1}:

\begin{thm}\label{poly1}
For every tree $T$, there is a polynomial $f(t)$ such that for every integer $t\ge 1$, if $G$ has no induced subgraph isomorphic to $T$ and 
no subgraph isomorphic to $K_{t,t}$, then $G$ has average degree at most $f(t)$.
\end{thm}

This is an improvement of a result of Kierstead and Penrice~\cite{kierstead}, who proved that there is a function $f(t)$ as in 
\ref{poly1}, not necessarily a polynomial; and that in turn was an improvement of a theorem of R\"odl (see \cite{gst,kierstead,kr}), 
who proved the same with average degree replaced by chromatic number.
Is there any hope for a comparable strengthening of \ref{mainthm}?

\ref{poly1} assumes that $G$ does not contain $K_{t,t}$ as a subgraph.
In \ref{mainthm} we replace this by the much weaker
hypothesis that the clique number of $G$ is bounded, although in compensation we must weaken the conclusion, replacing the bound on average degree with a bound on chromatic number.
This change is necessary: $K_{n,n}$ has large minimal degree, but no $K_3$ and no path-induced copy of $P_4$.  

But we 
could still ask for a {\em polynomial} bound on chromatic number.  Indeed, it is possible that the Gy\'arf\'as-Sumner conjecture holds 
with polynomial bounds (in other words, \ref{GSconj} with a bound on chromatic number that is polynomial in $t$).  This has recently 
been shown for a few trees, including every tree that does not contain the five-vertex path as an induced subgraph \cite{poly3}.  
However, the five-vertex path appears to be a sticking point.  Here, the best current bound is slightly superpolynomial 
(see \cite{poly4} for this and for related discussion):

\begin{thm}\label{p5}
If $G$ does not contain the five-vertex path $P_5$ as an induced subgraph, and has clique number $t$, then $\chi_G\le t^{\log_2 t}$.
\end{thm}
If $P$ is a path and $r$ is one end of $P$ then a graph $G$ contains a path-induced copy of $(P,r)$ if and only if it contains an induced copy of $P$.  Thus obtaining polynomial bounds in \ref{mainthm} even for $P_5$ would also require an improvement of \ref{p5}.

\end{document}